\documentclass[english]{smfart}
\usepackage{amssymb}
\usepackage[francais,english]{babel}
\usepackage{smfthm}
\usepackage[all]{xy}
\usepackage{bull}
\theoremstyle{plain}
\author{Luc Menichi}
\address{Universit\'e d'Angers, Facult\'e des Sciences\\
2 Boulevard Lavoisier\\49045 Angers, FRANCE,
Telephone: 02-41-73-50-25, Fax: 02-41-73-54-54}
\thanks{Topologie/Topology}
\email{Luc.Menichi@univ-angers.fr}
\title{P-th powers in mod p cohomology of fibers}
\alttitle{Les p-i\`emes puissances dans la cohomologie modulo p de fibres}
\begin{document}
\frontmatter
\begin{abstract}
Let $F\hookrightarrow E\twoheadrightarrow B$ be a fibration
whose base space $B$ is a finite simply-connected
CW-complex of dimension $\leq p$ and whose total space $E$ is a
path-connected CW-complex of dimension $\leq p-1$.
If $\alpha\in H^{+}(F;\mathbb{F}_p)$ then $\alpha ^{p}=0$.
\end{abstract}

\begin{altabstract}
Consid\'erons une fibration $F\hookrightarrow E\twoheadrightarrow B$ dont la
base $B$ est un CW-complexe fini simplement connexe de dimension $\leq
p$,
et dont l'espace total $E$ est un CW-complexe fini connexe par arcs de
dimension $\leq p-1$.
Si $\alpha\in H^{+}(F;\mathbb{F}_p)$ alors $\alpha ^{p}=0$.
\end{altabstract}
\maketitle

\mainmatter
We work over the prime field $\mathbb{F}_p$ with $p$ an odd or even prime.
The homology and cohomology of spaces are considered with
coefficients in $\mathbb{F}_p$.

In~\cite{AnickD:Hopah}, Anick proved using algebraic models:
\begin{enonce*}{Theorem~\cite[9.1]{AnickD:Hopah}}
Let $r$ be a non-negative integer.
Let $B$ be a simply-connected space with a finite type homology concentrated in degrees $i\in [r+1,rp]$.
Then all $p$-th powers vanish in $H^{+}(\Omega B)$.
\end{enonce*}

This result was suggested by McGibbon and Wilkerson~\cite[p. 699]{McGibbon-Wilkerson:loosfclp}.
The aim of this note is to give two different generalisations of Anick Theorem:
Theorem A and Theorem C below.

The first one, whose proof is inspired by a result of
Lannes and Schwartz~\cite[Prop 0.6]{Lannes-Schwartz:conSS},
uses the (vertical) Steenrod operations
in the Eilenberg-Moore spectral sequence:
\begin{enonce*}{Theorem A}
Let $r$ and $k$ be two non-negative integers.
Consider a fiber product of spaces
$$\xymatrix{
E\times_B X\ar[r]\ar@{>>}[d] &E\ar@{>>}[d]^{\pi}\\
X\ar[r] &B
}$$
where
\begin{itemize}
\item $\pi$ is a Serre fibration and
\item $H^{*}(E)$, $H^{*}(X)$ and $H^{*}(B)$ are of finite type.
\end{itemize}
If $B$ is simply-connected with homology concentrated in
degrees $i\in[r+1,rp^{k}]$
and the product space $E\times X$ is path connected
with homology $H_*(E\times X)$ concentrated in degrees $i\in
[r,rp^{k}-1]$,
then all $p^{k}$-th powers vanish in $H^{+}(E\times_B X)$.
\end{enonce*}
\begin{proof}
We suppose that $p$ is an odd prime. The case $p=2$ is similar.
Let $\mathcal{A}$ denote the mod $p$ Steenrod algebra.
The degree of an element $\alpha$ is denoted $\vert\alpha\vert$.
Recall from~\cite{RectorD:SteoEMss,SmithL:KuntEMss,SchwartzL:Unsmua}, that
the Eilenberg-Moore spectral sequence is a strongly convergent second
quadrant cohomological spectral sequence of $\mathcal{A}$-modules
$$E^{-s,*}_2\cong\mbox{Tor}^{-s,*}_{H^{*}(B)}\left(H^{*}(E),H^{*}(X)\right)
\Rightarrow H^{*}(E\times_B X).$$
More precisely, there exists a convergent filtration of
$\mathcal{A}$-modules on $H^{*}(E\times_B X)$:

$$H^{*}(E\times_B X)\supset\dots F_s\supset F_{s-1}\dots
F_1\supset F_0\supset F_{-1}=\{0\}$$
such that $\Sigma^{-s}F_s/F_{s-1}\cong E^{-s,*}_{\infty}, s\geq 0$.
Here $\Sigma^{-s}$ denotes the $s$-th desuspension of an
$\mathcal{A}$-module.

Let $\alpha\in F_s$ such that the class $[\alpha ]\in F_s/F_{s-1}$ is
non zero. We want to prove that $\alpha^{p^{k}}=0$.
As an $\mathcal{A}$-module,
$\mbox{Tor}^{-s,*}_{H^{*}(B)}\left(H^{*}(E),H^{*}(X)\right)$
is the $s$-th homology group of a
complex
of $\mathcal{A}$-modules, namely the Bar construction, whose
$s$-th term is $H^{*}(E)\otimes H^{+}(B)^{\otimes s}\otimes H^{*}(X)$.

The element $\Sigma^{-s}[\alpha ]\in E^{-s,*}_{\infty}$ is
represented by a cycle of the form $e[b_1|\dots|b_s]x$,
where $e\in H^{*}(E)$,
$(b_i)_{1\leq i\leq s}\in H^{+}(B)$ and $x\in H^{*}(X)$.
So $rs\leq\vert\alpha\vert\leq (rp^{k}-1)(s+1)$.

Case 1: $\vert e\vert+\vert x\vert\geq r$. Then $\vert\alpha\vert\geq r(s+1)$.
Therefore, by a degree argument, the element $\alpha^{p^{k}}$ of $F_s$ is
zero.

Case 2: $e=x=1$. Since the Cartan formula applies,
$\Sigma^{-s}[\alpha^{p}]=\mbox{P}^{\vert\alpha\vert/2}\Sigma^{-s}[\alpha]$
is represented by the element of the Bar construction,
$$\sum_{i_1+\dots+i_s=\vert\alpha\vert/2}[\mbox{P}^{i_1}b_1|\dots|\mbox{P}^{i_s}b_s].$$
So $\Sigma^{-s}[\alpha^{p^{k}}]$ is zero for degree reasons.
Therefore $\alpha^{p^{k}}$ belongs to $F_{s-1}$ which is concentrated
in degrees $\leq (rp^{k}-1)s$, thus $\alpha^{p^{k}}=0$.\end{proof}

\vspace{3mm}
Let $B$ be a space. The {\it free loop space} on $B$, denoted $B^{S^{1}}$,
is the set of continuous (unpointed) maps from the circle $S^{1}$ to $B$.
It can be defined as a fibre product:
$$\xymatrix{
B^{S^{1}}\ar@{^{(}->}[r]\ar@{>>}[d]_{ev} &B^{[0,1]}\ar@{>>}[d]^{\pi}\\
B\ar[r]_{\Delta} &B\times B
}$$
In this particular case, we can improve Theorem A.
\begin{enonce*}{Theorem B}
Let $r$ and $k$ be two non-negative integers.
If $B$ is a simply-connected space with finite type homology
concentrated in degrees $i\in [r+1,rp^{k}]$ then all
$p^{k}$-th powers vanish in $H^{+}(B^{S^{1}})$.
\end{enonce*}
\begin{proof} The Eilenberg-Moore spectral sequence for the previous fiber product
satisfies
$$E^{-s,*}_2\cong\mbox{HH}_s\left(H^{*}(B)\right)
\Rightarrow H^{*}(B^{S^{1}}).$$
Here $\mbox{HH}_*$ denotes the Hochschild homology.
As an $\mathcal{A}$-module,
$\mbox{HH}_s\left(H^{*}(B)\right)$
is the $s$-th homology group of a
complex
of $\mathcal{A}$-modules, namely the Hochschild complex, whose
$s$-th term is $H^{*}(B)\otimes H^{+}(B)^{\otimes s}$.

The same arguments as in the proof of Theorem A allow us to conclude except in Case 2
for $s=1$.
If $\alpha\in F_1\subset H^{*}(B^{S^{1}})$, we can only affirm that
$\alpha^{p^{k}}\in F_0$.
The evaluation map $ev:B^{S^{1}}\twoheadrightarrow B$ admits a section $\sigma$.
So $H^{*}(ev):H^{*}(B)\hookrightarrow (B^{S^{1}})$ admits $H^{*}(\sigma)$
as retract.
The edge homomorphism
$$H^{*}(B)=E^{0,*}_2\twoheadrightarrow E^{0,*}_3\dots \twoheadrightarrow E^{0,*}_\infty=F_0\subset H^{*}(B^{S^{1}})$$
correspond to $H^{*}(ev)$.
Since $\alpha^{p^{k}}\in F_0=H^{*}(B)$, $\alpha^{p^{k}}=\left[ H^{*}(\sigma)(\alpha)\right]^{p^{k}}$.
For degree reason, all $p^{k}$-th powers are zero in $H^{+}(B)$. So $\alpha^{p^{k}}=0$.\end{proof}

\vspace{3mm}
In~\cite{Felix-Halperin-Thomas:Serssmf},  F\'elix, Halperin and Thomas give a slightly
more complicated proof of Anick Theorem. Their proof uses the vertical and horizontal
Steenrod operations in the Serre spectral sequence:
\begin{enonce*}{Theorem~\cite[2.9(i)]{Felix-Halperin-Thomas:Serssmf}}
Let $r$ and $k$ be two non-negative integers.
If $B$ is a simply-connected space with a finite type homology
  concentrated in degrees $i\in [r+1,rp^{k}]$ then all $p^{k}$-th powers
  vanish in $H^{+}(\Omega B)$.
\end{enonce*}
This result generalizes in:
\begin{enonce*}{Theorem C \rm{(Compare with~\cite[10.8]{MenichiL:cohaf})}} 
Let $r$ and $k$ be two non-negative integers.
Let $F\buildrel{j}\over\hookrightarrow E\buildrel{\pi}\over\twoheadrightarrow B$
be a Serre fibration with $E$ path connected.
If $B$ is a simply-connected space with finite type homology
concentrated in degrees $i\in [r+1,rp^{k}]$ then, for any $\alpha\in H^{*}(F)$,
$\alpha^{p^{k}}\in\mbox{Im }H^{*}(j)$.
\end{enonce*}
\begin{proof}
The proof follows the lines of~\cite[2.9]{Felix-Halperin-Thomas:Serssmf}.
Since $H^{\leq r}(B)=0$, $\alpha\in E^{0,*}_2$ survives till $E^{0,*}_{r+1}$.
Therefore by a Theorem of Araki~\cite{ArakiS:Stepss} and V\'azquez~\cite{VasquezG:Stesss}
(See also~\cite[Prop 2.5 Case 2]{SawkaJ:oddpSoSS}), $\alpha^{p^{k}}\in E^{0,*}_2$
survives till $E^{0,*}_{rp^{k}+1}$.
Since $H^{>rp^{k}}(B)=0$, $$E^{0,*}_{rp^{k}+1}=E^{0,*}_\infty=\mbox{Im }H^{*}(j).$$\end{proof}

In order to see that the hypothesis in F\'elix-Halperin-Thomas Theorem
(and in Theorem B) cannot be improved,
consider $B=\Sigma\mathbb{CP} ^{p^{k}}$, the suspension of the $p^{k}$-dimension complex projective space.

Observe also that in Theorem C, $\alpha^{p^{k}}$ is not zero in general.
Indeed, take $\pi$ to be the fibration associated
to the suspension of the Hopf map from $S^{2p^{k}-1}$ to $\mathbb{CP} ^{p^{k}}
$~\cite[Remark 9.9]{MenichiL:cohaf}.

Finally, remark that the following question of
McGibbon and Wilkerson remains unsolved.
\begin{enonce*}{Question~\cite[p. 699]{McGibbon-Wilkerson:loosfclp} \rm{(See also~\cite[section 9.Question 3]{McGibbonC:pham})}}
Let $B$ be a finite simply-connected CW-complex and $p$ a prime large enough.
Do all the Steenrod operations act trivially on $H^{*}(\Omega B)$?
\end{enonce*}

I wish to thank Professor Katsuhiko Kuribayashi, for his precious help
with the Eilenberg-Moore spectral sequence.
\backmatter
\bibliography{Bibliographie}
\bibliographystyle{smfplain}
\end{document}